\documentclass{elsart}

\usepackage{epsfig}

\usepackage{amsmath,amsfonts}

\usepackage{amssymb}

\usepackage{rotating}

\usepackage{subfigure}

\newcommand{\fx} {\mathbf{x}}

\newcommand{\fXi} {\mathbf{\Xi}}
\newcommand{\fpsi} {\mathbf{\psi}}
\newcommand{\fa} {\mathbf{\alpha}}

\begin{document}

\begin{frontmatter}

\title{A maximum likelihood algorithm for the estimation and renormalization of exponential densities}

\author{Panagiotis Stinis}
\address{Department of Mathematics, 
Lawrence Berkeley National Laboratory, CA 94720 USA.}
\ead{stinis@math.lbl.gov}

\begin{abstract}
We present an algorithm based on maximum likelihood for the estimation and renormalization (marginalization) of exponential densities. The moment-matching problem resulting from the maximization of the likelihood is solved as an optimization problem using the Levenberg-Marquardt algorithm. In the case of renormalization, the moments needed to set up the moment-matching problem are evaluated using Swendsen's renormalization method. We focus on the renormalization version of the algorithm, where we demonstrate its use by computing the critical temperature of the two-dimensional Ising model. Possible applications of the algorithm are discussed.
\end{abstract}


\end{frontmatter}

\section{Introduction}
There are many problems in science and engineering that involve multiple 
time-scales or multiple spatial scales or both. However, due to computer limitations or due to the macroscopic nature of the quantities of interest,  
one wants to integrate out many of the variables involved in the system. 
This procedure involves two steps: i) one needs to specify the density of 
variables for the full system and ii) integrate out the unwanted variables 
(this amounts to the computation of the marginal probability density for the 
remaining variables). Both problems are of equal importance. One needs a 
good approximation for the density of {\it all} the variables involved in the system. 
But even if this is available, the computation of marginal probabilities can be costly. 
The present paper addresses both problems, focusing more on the second. Even 
though the algorithm is, in principle, equally applicable to the first problem,  we will 
present a detailed study of the behavior of the algorithm for that problem elsewhere. Interest in the 
evaluation of marginal probability densities is present in many different settings ranging 
from graphical models (the inference problem) \cite{jordan1} to reductions of systems 
of differential equations \cite{chorin5,chorin1}.

The algorithm is based on maximum likelihood estimation. We should mention here the algorithm of Geyer and Thompson \cite{geyer} for the estimation of the parameters of an unknown exponential density, which is also based on maximum likelihood but whose approach is different from ours (see also \cite{besag}). To the best of our knowledge, the algorithm that we present here for renormalization is the first one to be based solely on maximum likelihood (see also \cite{chorin5} for a different approach based on conditional expectations). As will be explained more below, the numerical implementation of the algorithm requires the solution of an optimization problem. The solution of the optimization problem through the Levenberg-Marquardt algorithm is efficient and robust. It, also, avoids some problems associated with more conventional methods like steepest descent and Newton's. The advantage of the Levenberg-Marquardt algorithm becomes especially crucial in the case of renormalization, where one needs to determine very accurately the parameters of the renormalized (marginal) exponential density. For the example we present (the $2D$ Ising model of spins), this allows the accurate determination of the critical temperature.

Exponential densities are, partly due to their nice mathematical features, widely used in the modeling of densities of systems of interacting variables in different contexts, ranging from Hamiltonian systems to image processing and bioinformatics (see \cite{jordan2} and references therein). As a result, there is an increased interest in algorithms for estimating and manipulating such densities numerically. What we offer here is a general algorithm that allows the estimation of parameters of exponential densities. In addition to estimating the parameters of an exponential density, and this is the main focus in the present work, it allows the renormalization of a known exponential density. Renormalization amounts to calculation of marginal densities in a way that the functional form of the density is retained. Since we want to retain the mathematical structure of the density, the marginal density that we compute will not always be exact. There is a trade-off between efficiency and accuracy in such calculations that is well known in the context of real-space renormalization in statistical physics \cite{binney,swendsen1}. For the test case  presented here (the $2D$ Ising model of interacting spin variables), this does not turn out to be harmful.

Suppose we are given a number of independent samples from a density and we try to fit an exponential density to these samples by maximizing their likelihood. Whether we are looking for the parameters of an unknown exponential density or the renormalized parameters of a known exponential density, maximum likelihood estimation of these parameters leads to a moment-matching problem. In other words, we want to determine the parameters of an exponential density so that a finite number of its moments match the moments computed from the given samples. In the present context, the word moment stands for the empirical average (expectation value) of a, not necessarily polynomial, function of the random variables. For the problem of estimating the parameters of an unknown exponential density, the moments needed to set up the moment-matching problem are computed by using the samples of the unknown density to which we are trying to fit the exponential density. For the problem of estimating the renormalized parameters of a known exponential density, the moments needed to set up the matching problem are computed using Swendsen's renormalization method \cite{swendsen1}. In both cases, the equations that define the moment-matching problem contain, in general, nonlinear functions of the parameters to be estimated. We solve the matching problem as an optimization problem using the Levenberg-Marquardt algorithm (see e.g. \cite{bishop}). The Levenberg-Marquardt algorithm is a combination of Newton's method and the method of steepest descent. The reason for using this algorithm is that it allows flexibility in the initial guess of the parameters. Use of the steepest descent algorithm alone can lead to slow convergence, while use of Newton's method alone can lead to divergence, because most likely the initial guess of the parameters is not close to their true values. 

The paper is organized as follows. In Section \ref{algo} we present the algorithm for the estimation of the parameters of an unknown exponential density. We also present the modifications needed to renormalize a known exponential density. In Section \ref{num} the algorithm is used, first to estimate, and then to renormalize  the parameters of the two-dimensional Ising model of $\pm 1$ spins. The estimates of the renormalized parameters are used to locate the critical temperature of the model. In the final section, we discuss possible applications of the algorithm. 

\section{The algorithm}{\label{algo}}
In this section we present an algorithm that allows the estimation of parameters of exponential densities. We examine two cases. In the first case we are given a number of samples from a multivariate density and we estimate the parameters for an exponential representation of the density. In the second case, we are given the parameters of a multivariate exponential density and we are presented with the problem of computing marginal probabilities so that the form of the density remains the same.

\subsection{Estimation of parameters for an unknown exponential density}{\label{unkno}}
We begin our presentation with a few facts about families of exponential densities and convex analysis (see \cite{brown,rocka,jordan2}). Let $\mathbf{x}=(x_1,\ldots,x_n)$ be an $n$-dimensional random vector taking values in $\mathbf{\Xi}^n \subseteq \mathbb{R}^n.$ The set $\mathbf{\Xi}^n=\Xi_1 \times \Xi_2 \times \ldots \Xi_n,$ where $x_1 \in \Xi_1,\ldots , x_n \in \Xi_n.$ Also, let $\psi_k(\mathbf{x}), k=1,\ldots,l$ be a collection of functions of $\mathbf{x}.$ The functions $\psi_k$ are known as potentials or sufficient statistics. Let $\psi=(\psi_1,\ldots,\psi_l)$ be the vector of potential functions. Associated with the vector $\psi$ is a vector $\mathbf{\alpha}=(\alpha_1,\ldots,\alpha_l)$ whose elements are called canonical or exponential parameters. The exponential family associated with $\mathbf{\psi}$ is the collection of density functions (parametrized by $\mathbf{\alpha}$) of the form
$$p(\mathbf{x},\mathbf{\alpha})=\frac{\exp(-\langle \mathbf{\alpha},\mathbf{\psi}(\mathbf{x})\rangle)}{Z(\mathbf{\alpha})},$$
where $\langle \fa,\mathbf{\psi}(\fx) \rangle=\sum_{k=1}^l \alpha_k \psi_k(\fx)$  and $Z(\mathbf{\alpha})=\int_{\mathbf{\Xi}^n} \exp(-\langle \mathbf{\alpha},\mathbf{\psi}(\mathbf{x})\rangle)d\mathbf{x}.$ The exponential family is defined only for the set 
$$ \mathbf{A}= \{ \mathbf{\alpha} \in \mathbb{R}^l | Z(\mathbf{\alpha}) < \infty \}.$$
If $\mathbf{A}$ is open, the exponential family is called regular. We will restrict our attention to regular families, so that in all the theorems stated below the assumption of regularity will be implied. Usually the exponential density is defined without the minus sign, but we incorporate it in anticipation of the case of Hamiltonian systems like the Ising model that we examine in Section \ref{num}. In that case we have $H(\mathbf{x})=\langle \mathbf{\alpha},\mathbf{\psi}(\mathbf{x})\rangle,$ where $H(\mathbf{x})$ is the Hamiltonian of the system and thus, the exponential density is the Boltzmann density. In the case where the variables $x_i, i=1,\ldots,n$ take only a discrete number of values, the integration is replaced by summation over the possible values for each of the $x_i.$ If the functions $\psi_k(\mathbf{x})$ are linearly independent,  the representation is called minimal. Otherwise, it is called overcomplete. The distinction between minimal and overcomplete representations will be used later when we formulate the moment-matching problem.

Suppose that we are given a collection of $N$ independent samples of an $n$-dimensional random vector $\mathbf{x}$. In general, we do not know which density the samples are drawn from. There are many examples in practical applications where the random vector comes from a exponential density (see \cite{jordan2}). However, even if we do not know that the samples are drawn from an exponential density, we can try to fit an exponential density to the samples. This will become clearer when we formulate the moment-matching problem and exploit some properties of the exponential densities. The basic idea behind the algorithm we present here is to estimate the unknown parameter vector $\mathbf{\alpha}$ by maximizing the likelihood function of the samples. For a collection of $N$ independent samples of the random vector $\mathbf{x}$, the likelihood function $L$ is defined as (see e.g. \cite{lehman}) 
$$ L=\prod_{j=1}^{N} p(\mathbf{x}_j,\alpha),$$ 
where  $p(\mathbf{x}_j,\alpha)$ is the unknown exponential density whose parameters $\mathbf{\alpha}$ we wish to determine. We associate a potential function $\psi_k, k=1,\ldots,l$ with every parameter $\alpha_k.$ Maximization of $L$ with respect to the parameters $\alpha_k$ produces an estimate $\bar{\mathbf{\alpha}}$ for $\mathbf{\alpha}.$  Under suitable regularity conditions, the sequence of estimates $\bar{\mathbf{\alpha}}$ for increasing values of $N$ is asymptotically efficient and tends, with probability one, to a local maximum in parameter space. From now on we will use the notation $\mathbf{\alpha}$ instead of $\bar{\mathbf{\alpha}}$ to denote the maximum likelihood estimate of the parameters keeping in mind that this is only an estimate of the parameters. In addition, we will be working with the logarithm of the likelihood $\log L$, since it does not alter the position of the maximum and also leads to formulas that are more easily manipulated. Differentiation of $\log L$ with respect to the $\alpha_k$ results in 
\begin{equation}
\label{mom1}
E_{\mathbf{\alpha}}[\psi_k(\mathbf{x})]=\frac{1}{N} \sum_{j=1}^{N} \psi_k(\mathbf{x}_j), \quad k=1,\ldots,l
\end{equation}
where
$$E_{\mathbf{\alpha}}[\psi_k(\mathbf{x})]= \frac{\int_{\mathbf{\Xi}^n} \psi_k(\fx) \exp(-\langle \mathbf{\alpha},\mathbf{\psi}(\mathbf{x})\rangle) d\mathbf{x}}{\int_{\mathbf{\Xi}^n} \exp(-\langle \mathbf{\alpha},\mathbf{\psi}(\mathbf{x})\rangle) d\mathbf{x}} $$ 
is the expectation value of the function $\psi_k$ with respect to the density $p(\mathbf{x},\mathbf{\alpha}).$ The right side of (\ref{mom1}) is the average (moment) of the function $\psi_k$ as calculated from the given samples. The $l$ equations in (\ref{mom1}) define the {\it moment-matching} problem. What we want to do is to estimate the parameters $\mathbf{\alpha}$ so that the conditions in (\ref{mom1}) are satisfied. The question is whether such a problem has a solution, and if it does whether it is unique. To answer this we resort to the theory of exponential densities and convex analysis.

First we should note that the moments of the potential functions define an alternative parametrization of the exponential family. This is known as mean parametrization. In fact, let $\mathbf{\mu} \in \mathbb{R}^l$ be the vector of moments of the potential functions. Also, define the set $M$ as
$$M= \{ \mathbf{\mu} \in \mathbb{R}^l | \exists p(\cdot) : \int_{\mathbf{\Xi}^n} \mathbf{\psi}(\mathbf{x}) p(\mathbf{x}) d\mathbf{x}=\mathbf{\mu} \}.$$
Note that $M$ is a convex set and that in the definition of $M$ we do not restrict the density $p(\cdot)$ to the exponential family. The density $p(\cdot)$ is {\it any} density that realizes $\mathbf{\mu}.$ Typically, the exponential family $\{ p(\mathbf{x},\mathbf{\alpha}) | \mathbf{\alpha} \in \mathbf{A} \}$ is only a strict subset of all possible densities. Since we are interested in exponential densities, we have to find the relation between the set $\mathbf{A}$ of admissible parameter vectors and the set $M.$ This will allow us to answer the question whether the moment-matching problem admits  an exponential density as a solution.

For a given vector of potentials $\mathbf{\psi},$ define the mapping $\Lambda : \mathbf{A} \rightarrow M$ as
$$\Lambda(\mathbf{\alpha})=E_{\mathbf{\alpha}}[\mathbf{\psi}(\mathbf{x})]=\int_{\mathbf{\Xi}^n} \mathbf{\psi}(\mathbf{x}) p(\mathbf{x},\fa) d\mathbf{x}$$
Whether there exists a parameter vector $\fa$ satisfying (\ref{mom1}) depends on the properties of the mapping $\Lambda.$ In particular, it depends on: i) if $\Lambda$ is one-to-one and hence invertible on its image and ii) what is the image of $\mathbf{A}$ under $\Lambda.$ There are two theorems that characterize the properties of $\Lambda$ (see \cite{jordan2}). For the first theorem we need the distinction between a minimal and an overcomplete representation.
\begin{thm}
\label{t1}
The mapping $\Lambda$ is one-to-one if and only if the exponential representation is minimal. 
\end{thm}

For the second theorem we need the definition of the relative interior of a set. The relative interior of a set is the interior taken with respect to its affine hull. A key result from convex analysis states that for any non-empty convex set its relative interior is non-empty.
\begin{thm}
\label{t2}
The mapping $\Lambda$ is onto the relative interior of $M$. 
\end{thm}
Theorem \ref{t2} in conjunction with Theorem \ref{t1} guarantees, for minimal exponential representations, the existence of a {\it unique} parameter vector for each point in the relative interior of $M.$ Of course, there is the question of what happens for points in the closure of $M$ that are not in the relative interior. To answer that we use one more result from convex analysis.
\begin{thm}
\label{t3}
Let $M$ be a convex set in $\mathbb{R}^l$. Let $x$ a point in the relative interior of $M$ and $y$ a point in the closure of $M.$ Then $\lambda x +(1-\lambda)y$ belongs to the relative interior of $M$ for $0 < \lambda  \leq 1.$
\end{thm}
As we have mentioned before, the exponential family typically describes only a strict subset of all possible densities that give rise to the set $M.$ However, Theorems \ref{t1}-\ref{t3} tell us that this is enough for the moment-matching problem. In the case of an overcomplete representation, there is no longer a one-to-one correspondence between $\mathbf{A}$ and $\Lambda(\mathbf{A}).$ But this is not a problem. For an overcomplete representation, the solution for the moment-matching problem is no longer unique but it exists. Any of the solutions are equally admissible, since all of them reproduce the same moments.

Now that we have defined the moment-matching problem we have to find a way to actually estimate the parameter vector $\fa.$ The equations (\ref{mom1}) contain, in general, nonlinear functions of the parameters. Moreover, except for very special cases, these nonlinear functions are unknown or very difficult to manipulate analytically. Thus, we have to tackle the problem of estimating the parameter vector numerically. We can define the $l$-dimensional vector  $\mathbf{f}(\fa)=(f_1(\fa),f_2(\fa),\ldots,f_l(\fa))$ as
\begin{equation}
\label{mom2}
f_k(\fa)=E_{\mathbf{\alpha}}[\psi_k(\mathbf{x})]-\frac{1}{N} \sum_{j=1}^{N} \psi_k(\mathbf{x}_j), \quad k=1,\ldots,l
\end{equation}
The moment-matching problem amounts to solving the system of (nonlinear) equations $f_k(\fa)=0, k=1,\ldots,l.$ Two popular candidates to perform such a task are the method of steepest descent and Newton's method. However, both have their drawbacks. The method of steepest descent converges but can have very slow convergence, while Newton's method converges quadratically but it diverges if the initial guess of the solution is not good. We choose to solve the moment-matching problem as an optimization problem using the Levenberg-Marquardt (LM) algorithm (see e.g. \cite{bishop}). This is a powerful iterative optimization algorithm that combines the advantages of the method of steepest descent and Newton's method. First, let us write the moment-matching problem as an optimization problem. Define the error function $\epsilon(\fa)$ as
$$\epsilon(\fa)=\frac{1}{2}\sum_{k=1}^l \epsilon_k^2=\frac{1}{2}\sum_{k=1}^l f_k^2(\fa),$$
where $\epsilon_k=f_k(\alpha).$
The problem of minimizing $\epsilon(\fa)$ is equivalent to solving the system of equations $f_k(\fa)=0, k=1,\ldots,l$ i.e. the zeros of $\epsilon$ are solutions of the system $f_k(\fa)=0, k=1,\ldots,l$ and vice versa. The LM algorithm uses a positive parameter $\lambda$ to control convergence and the updates of the parameters at step $m+1$ are calculated through the formula
\begin{equation}
\label{lm1}
\alpha_k^{m+1}={\alpha_k^m} -[J^T J+ \lambda  diag (J^T J) ]^{-1} J^T  \mathbf{f}(\fa^m),
\end{equation}
where $J=\frac{\partial{f_i}}{\partial {\alpha_j}}|_{\fa=\fa^m}, i,j=1,\ldots,l$ is the Jacobian of $\mathbf{f}(\fa^m)$ and $J^T$ its transpose. The matrix $diag (J^T J)$ is a diagonal matrix whose diagonal elements are the diagonal elements of $ (J^T J).$ In the literature, the name Levenberg-Marquardt is also used to denote the algorithm in (\ref{lm1}) with $diag (J^T J)$ replaced by the unit matrix $I.$ For the case where we use the unit matrix $I$ instead of $diag (J^T J),$ it is straightforward to see the connection with the methods of steepest descent and Newton's. For $\lambda=0$ the algorithm reduces to Newton's method, while for very large $\lambda$ we recover the steepest descent method. The modification (due to Marquardt) of using $diag (J^T J)$ becomes important in the case where $\lambda$ is large. In this case if we only used the unit matrix $I$ almost  all information coming from $ (J^T J)$ is lost. On the other hand, since $ (J^T J)$ provides information about the curvature of $\epsilon$, use of the matrix $ diag (J^T J)$ allows us to incorporate information about the curvature even in cases with large $\lambda.$ In the numerical simulations we used both forms of the algorithm. The form in (\ref{lm1}) gave superior results.

We have to prescribe a way of computing the Jacobian $J(\fa^m)$. The element $J_{ij}$ of the  Jacobian is given by 
$$J_{ij}(\fa^m)=- (E_{\fa^m}[\psi_i(\fx)\psi_j(\fx)]-E_{\fa^m}[\psi_i(\fx)]E_{\fa^m}[\psi_j(\fx)] )$$
for $i,j=1,\ldots,l$ (note that the Jacobian is symmetric) So, all the quantities involved in equation (\ref{lm1}) can be expressed as expectation values with respect to the $m$-th step parameter estimate $\fa^m.$ We compute these expectation values using the Metropolis Monte Carlo algorithm. This can make the algorithm expensive since the density  has to be sampled at each step. However, the cost of the algorithm can be reduced by parallelization of the Monte Carlo sampling procedure. Also, note that if one uses as potential functions (non-orthogonal) polynomials, the condition number of the Jacobian matrix grows fast with the order of polynomials. An ill-conditioned Jacobian can lead to catastrophic errors in the evaluation of the parameter vector. This is especially crucial in the renormalization version of the algorithm, where one wants to use the parameter vector to look for possible phase transitions and their associated critical properties. This point will become clearer in Section \ref{num} where we present results for the $2D$ Ising model of spins.

We conclude this section with some comments about the value of $\lambda$ and the convergence criterion used to stop the iterative process. Since $\lambda$ acts as a regulator between the steepest descent and Newton aspects of the algorithm, its value should be determined in a way that brings out the advantages of the two methods. The starting value of $\lambda$ was chosen to be 1. When we detected a streak of a few error-decreasing steps, the value of $\lambda$ was decreased by a factor of 10. On the other hand, if the error increased, we repeated the step with the value of $\lambda$ increased by a factor of 10. We also need to prescribe a convergence criterion. We used the relative error criterion, i.e. the algorithm is stopped if $|(\epsilon_{new}-\epsilon_{old})/\epsilon_{old} |< RTOL.$ The value of $RTOL$ is determined by the accuracy of the Monte Carlo sampling. The relative error criterion is not enough on its own because there is the possibility of a very small step while the algorithm is still far away from the minimum of the error. Such a small step would pass the relative error criterion and the algorithm would stop. To avoid such false convergence, we added an extra convergence check criterion. Whenever the relative error criterion was satisfied, we checked that the error value was acceptable under the absolute value criterion $\underset{1 \leq k \leq l}{max}\sqrt{(|\epsilon_k|/2)}/{|\mu_k|} < ATOL$ where $\mu_k=\frac{1}{N} \sum_{j=1}^{N} \psi_k(\mathbf{x}_j).$ Note that the absolute error criterion has to take into account the magnitude of the moment. Also, $ATOL$ cannot be smaller than the accuracy afforded by the Monte Carlo sampling. If the error did not pass the absolute error criterion, the algorithm was not stopped even if the relative criterion was satisfied. Finally, we added a constraint on the total number of iterations allowed and on the maximum value of $\lambda.$

\subsection{Renormalization (marginalization) of a known exponential density}{\label{kno}}

We present the necessary modifications to the above scheme in the case when we know the parameters of an exponential density and we want to compute marginal densities. Suppose that we know the parameter vector $\fa$ for a given vector $\fpsi(\fx)$ of potential functions for an $n$-dimensional random vector $\fx.$ The exponential density associated with these parameters and potentials is $\frac{\exp(-\langle \mathbf{\alpha},\mathbf{\psi}(\mathbf{x})\rangle)}{Z(\mathbf{\alpha})}.$ Suppose that we want to retain only the first $m$ variables. The vector $\fx$ can be written as $\fx=(\hat{\fx},\tilde{\fx}),$ where $\hat{\fx}$ is the vector of the first $m$ variables and $\tilde{\fx}$ is the vector of the remaining $n-m$ variables. Denote the exponent as $H(\mathbf{x})=\langle \mathbf{\alpha},\mathbf{\psi}(\mathbf{x})\rangle$ (this is just a notation, it does not imply that we are only considering Hamiltonian systems). The problem of finding a function 
$\hat{H}(\hat{\fx})$ such that 
$$\exp(-\hat{H}(\hat{\fx}))=\int_{\fXi^{n-m}} \exp(-H(\fx)) d\tilde{\fx}$$
is well-defined, at least for a vector $\fx$ of finite dimensionality. However, as we have already mentioned in the introduction, we not only want to compute exponential marginal densities, but do so while retaining the mathematical form of the exponent. This means that we want the marginal density's exponent $\hat{H}(\hat{\fx})$ to have the form $\hat{H}(\hat{\fx})=\langle \hat{\mathbf{\alpha}},\hat{\mathbf{\psi}}(\hat{\mathbf{x}})\rangle.$ The functions $\hat{\psi}_k(\hat{\fx}), \, k=1,\ldots,l$ have the same form as the functions $\psi_k(\fx),$ but are defined {\it only} over the $m$ variables. We can think of the representation of the exponent through the potential functions as an expansion of the functions $H$ and $\hat{H}$ in functions of $\fx$ and $\hat{\fx}$ respectively. The functions $H$ and $\hat{H}$ are expanded using the same form of potential functions (on different sets of variables). It is clear that the way we have chosen to represent the marginal density can result in this density being only approximate. The reason is that the marginal density may involve potential functions in addition, or different, from the potential functions that we have chosen. How much of a problem this will turn out to be, depends on our choice of potential functions.  If one has no physical intuition or prior knowledge about the form of the potential functions, one can use as potential functions polynomials in the reduced vector $\hat{\fx}.$ For example, if the random vector is continuous, one could use products of one-dimensional Legendre polynomials. Usually, one augments the vector of potential functions of the original density by adding components with the corresponding parameters set to zero. The additional parameters may or may not acquire nonzero values for the marginal density, depending on how well we choose them.

Thus, the problem we are addressing is whether we can find a parameter vector $\hat{\fa}$ such that 
\begin{equation}
\label{den1}
\exp(-\langle \hat{\mathbf{\alpha}},\hat{\mathbf{\psi}}(\hat{\fx})\rangle)=\int_{\fXi^{n-m}} \exp(-\langle \mathbf{\alpha},\mathbf{\psi}(\mathbf{x})\rangle) d\tilde{\fx}.
\end{equation}

If (\ref{den1}) holds, then the normalization constant $Z$ satisfies $Z(\hat{\fa})=Z(\fa),$ since
\begin{equation*}
\begin{split}
Z(\hat{\fa})=\int_{{\fXi}^m} \exp(-\langle \hat{\mathbf{\alpha}},\hat{\mathbf{\psi}}(\hat{\fx})\rangle) d\hat{\fx} &=  \int_{{\fXi}^m} \biggl [  \int_{{\fXi}^{n-m}} \exp(-\langle \hat{\mathbf{\alpha}},\hat{\mathbf{\psi}}({\fx})\rangle)d\tilde{\fx} \biggr ] d\hat{\fx} \\
&= \int_{{\fXi}^n} \exp(-\langle {\mathbf{\alpha}},\mathbf{\psi}({\fx})\rangle)d{\fx} =Z(\fa) \\
\end{split}
\end{equation*}
Equation (\ref{den1}) can be used again to compute the marginal density for a subset $\hat{\hat{\fx}}$ of $\hat{\fx}$ and so on. Also it can be applied to more general groupings of variables in $\fx,$ e.g. the renormalization scheme that we will use in Section \ref{num} (see also \cite{binney}).  

Now we turn to the problem of estimating the parameter vector $\hat{\fa}.$ Since the marginal density is also of exponential form, we can apply the algorithm of Section \ref{unkno}. The resulting moment-matching problem is 
\begin{equation}
\label{mom1r}
E_{\hat{\fa}}[\hat{\psi}_k(\hat{\fx})]=\frac{1}{N} \sum_{j=1}^{N} \hat{\psi}_k(\hat{\fx}_j), \quad k=1,\ldots,l
\end{equation}
where
$$E_{\hat{\fa}}[\hat{\psi}_k(\hat{\fx})]= \frac{\int_{\mathbf{\Xi}^m} \hat{\psi}_k(\hat{\fx}) \exp(-\langle \hat{\fa},\hat{\mathbf{\psi}}(\hat{\fx})\rangle) d\hat{\fx}}{\int_{\mathbf{\Xi}^m} \exp(-\langle \hat{\fa},\hat{\mathbf{\psi}}(\hat{\fx})\rangle) d\hat{\fx}} $$ 
is the expectation value of $\psi_k$ with respect to the density
$\exp(-\langle \hat{\fa},\hat{\fpsi}(\hat{\fx})\rangle)/Z(\hat{\fa}).$ If we look at the equations (\ref{mom1r}), we see that we need samples of the random vector $\hat{\fx}.$ This can be effected through Swendsen's observation \cite{swendsen1}, that the marginal density can be sampled without knowing its explicit form. Since we know the density for the vector $\fx,$ we can sample it using e.g. Metropolis Monte Carlo and obtain samples of the vector $\hat{\fx}$ by keeping only the first $m$ variables of each vector $\fx.$ After we obtain the samples of the vector $\hat{\fx}$ we can apply the rest of the algorithm as it is and estimate the parameter vector $\hat{\fa}.$ This procedure can be performed recursively (as is done e.g. in real-space renormalization \cite{binney}), and thus obtain a "parameter" flow which contains useful information about the system. This is done in the next section for the two-dimensional Ising model of spins.

Note that if we are interested only in producing samples of the reduced vector $\hat{\fx},$ Swendsen's method suffices. However, to use Swendsen's method we have to compute first samples from the $n$-dimensional vector $\fx$ and this can be very costly. What we offer here is a way of representing the marginal density for the reduced vector $\hat{\fx}$ with an analytical formula. All that we need to describe the marginal density is the parameter vector $\hat{\fa}.$ This is highly desirable if e.g. one wants to sample the marginal density in the future or compute conditional expectations with respect to the marginal density. Suppose that we want to compute the conditional expectation of a function $h(\hat{\fx})$ with respect to a subset of the (already) reduced vector $\hat{\fx}.$ Let the reduced vector be split up as $\hat{\fx}=(\hat{\fx}^1,\hat{\fx}^2)$ where $\hat{\fx}^1$ is $m_1$-dimensional and $\hat{\fx}^2$ is of dimension $m_2=m-m_1.$ Suppose that we want to compute the conditional expectation of $h(\hat{\fx})$ conditioned on $\hat{\fx}^1.$ Since we know the analytical expression for $p(\hat{\fx},\hat{\fa}),$ this amounts to the calculation (e.g. using Monte Carlo) of the quantity
\begin{equation}
\label{con1r}
\begin{split}
E_{\hat{\fa}}[h(\hat{\fx})|\hat{\fx}^1] & =\frac{\int_{{\fXi}^{m_2}} h(\hat{\fx}) \exp(-\langle \hat{\mathbf{\alpha}},\hat{\mathbf{\psi}}(\hat{\fx})\rangle) d\hat{\fx}^2}{\int_{{\fXi}^{m_2}} \exp(-\langle \hat{\mathbf{\alpha}},\hat{\mathbf{\psi}}(\hat{\fx})\rangle) d\hat{\fx}^2} \\
& =\frac{\int_{{\fXi}^{n-m+m_2}} h(\hat{\fx}) \exp(-\langle {\mathbf{\alpha}},\mathbf{\psi}({\fx})\rangle) d\tilde{\fx} d\hat{\fx}^2}{\int_{{\fXi}^{n-m+m_2}} \exp(-\langle {\mathbf{\alpha}},\mathbf{\psi}({\fx})\rangle) d\tilde{\fx} d\hat{\fx}^2}\\
\end{split}
\end{equation}
where the second equality follows from (\ref{den1}).
If we did not have an analytical formula we would have to sample the original $n$-dimensional vector fixing the values of $\hat{\fx}^1,$ i.e. use the second equality in (\ref{con1r}). This can be very costly when $n$ is large (as it happens usually in applications). On the other hand, the computation in the first equality in (\ref{con1r}) is much cheaper. This can be useful e.g. in the inference problem for a graph with cycles (in the context of graphical models), where one can compute the marginal densities on the different cliques (fully-connected clusters of nodes) and store the parameter vectors for further calculations within the individual cliques. In fact, the process can be parallelized, by assigning a processor to one (or a few) cliques. It would be interesting to see how this algorithm compares with the junction tree algorithm (see \cite{lauritzen}) whose complexity grows exponentially with the size of the maximal clique. The present algorithm can be applied directly on the clique tree (an acyclic graph whose nodes are formed by cliques of the original graph) without the need for triangulation. It is also interesting to see how the present algorithm compares with variational inference methods \cite{jordan2}.

From equations (\ref{mom1}) and (\ref{mom1r}), we see that the problem of maximum likelihood estimation of the parameter vector for both the full system and the reduced system is transformed into the same moment-matching problem, i.e. involving the same number and same form of potential functions. The only difference is in the dimensionality of the random vector involved.

\section{Numerical results for the $2D$ Ising spin model}{\label{num}}

In this section we apply the algorithm from Section \ref{algo} to the $2D$ Ising spin model. Since the system has a known exponential probability density, we can check the accuracy of the algorithm in computing this density starting from a different initial guess. Also, the Ising model exhibits a phase transition with known properties. We use the "renormalized" version of the algorithm (Section \ref{kno}) to compute the critical temperature.

Consider a square lattice of size $L$. We denote a lattice site as $I_k=(i_k,j_k)$ where $i_k,j_k$ are integers and $k=1,\ldots,L^2$ (we assume that the sites are listed in a convenient order). We associate with each lattice site a (spin) variable $x_k$ that can only take the values $\pm 1.$ The variables on the whole lattice will be denoted by $\fx=(x_1,\ldots,x_{L^2}).$ The lattice is periodic with period $L.$ The Ising model of interaction of the variables (see e.g. \cite{binney}) is defined through the Hamiltonian $$H=-\frac{1}{T} \sum_{<I,J>}x_I x_J$$ where $<>$ means summation only over nearest neighbors. The parameter $T$ is interpreted as the temperature in the context of statistical mechanics. The model's probability density function is $$\frac{1}{Z}\exp(-H),$$ where $Z$ is the normalization constant (partition function). The Ising model has been used extensively to study phase transitions (as a function of the temperature) in magnetic materials. According to Onsager's analytical solution, the average magnetization $m=E[\sum_{I=1}^{L^2} x_I/{L^2}]$ exhibits, in the limit of $L \rightarrow \infty,$ a transition from zero to non-zero values for temperatures $T \leq T_c=2.269.$

We can cast the probability density of the Ising model in the form of an exponential density as the ones described in Section \ref{algo}. To do that we will need to define groups of variables around a site, say $I=(i,j).$ These groups have nothing to do with the groupings of variables that we will use later for real-space renormalization. The groups are defined as follows:  group 1 contains only $x_I.$ Group 2 contains the variables whose distance from $I$ is 1 (the nearest neighbors), group 3 contains those variables whose distance is $\sqrt{2}$, group 4 contains the variables whose distance is $2$ etc. We use the members of the groups to define the corresponding collective variables $X_{I,k}$ as $$X_{I,k}=\frac{1}{n_k} \sum_{group \,  k} x_J,$$ where $n_k$ is the number of variables in group $k.$ These collective variables can be used to form translation-invariant polynomials in $\fx,$ for example,  $\sum_{I=1}^{L^2} X_{I,1} (X_{I,k})^p.$ Using this notation, the Hamiltonian for the Ising model can be written as $H=-\frac{2}{T}\sum_{I=1}^{L^2} X_{I,1} X_{I,2}=\frac{2}{T} \bigl ( -\sum_{I=1}^{L^2} X_{I,1} X_{I,2}\bigr )$ For the case of the Ising model, we will include a minus sign in the definition of all the potential functions. The thermodynamic properties are the same whether or not we include the minus sign \cite{binney}. In the numerical simulations we used the following 8 potential functions
\begin{equation}
\label{pots}
\begin{split}
\psi_k=&-\sum_J X_{J,1} X_{J,k+1}   \quad \text {for} \: k\,=\, 1,2,3,4,5 \\
\psi_{k+5}=&-\sum_J (X_{J,k+1})^4 \quad \text {for} \: k\,=\, 1,2 \\
\psi_8=&-\sum_J X_{J,2}^2 X_{J,3}^2,
\end{split}
\end{equation}
where the summation extends over the whole lattice (see also \cite{chorin5}). The Hamiltonian for the Ising model can be written as $H(\alpha,\fx)=\sum_{k=1}^{8} \alpha_k \psi_k (\fx),$ with $\alpha_1=2/T$ and $\alpha_l=0$ for $l=2,\ldots,8.$ The coefficients for the potential functions $\psi_2 \ldots \psi_8$ acquire nonzero values when we compute marginal densities for subsets of $\fx.$

\begin{figure}
\centering
\epsfig{file=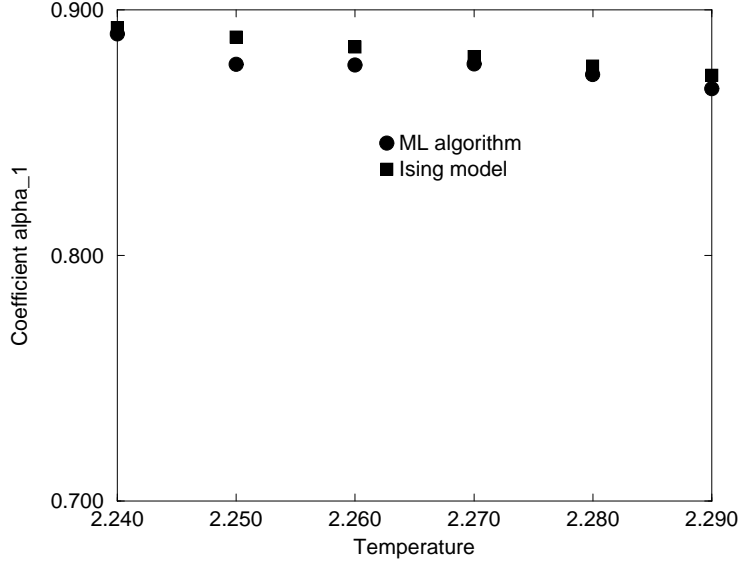,height=3in}
\caption{Comparison of the value of the parameter $\alpha_1$ as determined by the algorithm with its true value.}
\label{fig:unkno1}
\end{figure}

The algorithm we described in Section \ref{algo} can be used to estimate the parameters of an unknown exponential density or to renormalize a known density. As a test of the first aspect of the algorithm we compute the parameters for the Ising model. In other words, we prescribe some initial values for the parameters $\alpha_k, \, k=1,\ldots,8$ and we apply the algorithm. The parameters computed by the algorithm should converge to the true values $\alpha_1=2/T$ and $\alpha_l=0$ for $l=2,\ldots,8.$ We should note here, that the representation we are using is not minimal and thus, there are more than one admissible sets of parameters for the same moments. In the case of the Ising model, if we want to recover the values $\alpha_1=2/T$ and $\alpha_l=0$ for $l=2,\ldots,8$ we should start close to these values to ensure that we are in their basin of attraction. Of course, any set of parameters that is produced from a convergent run of the algorithm is equally acceptable. In Fig.(\ref{fig:unkno1}) we compare the value of the parameter $\alpha_1$ as determined by the algorithm with its true value $2/T$ for different temperatures. The initial guess was $1.9/T$ for $\alpha_1$ and zero for the rest of the parameters. The results shown are for $10^5$ Monte Carlo steps per spin (each step involves a sweeping of the whole lattice). The agreement for $\alpha_1$ shown in Fig.(\ref{fig:unkno1}) is within the accuracy of the Monte Carlo computation. 
The values of the rest of the parameters $\alpha_2,\ldots,\alpha_8,$ remain zero to within the accuracy afforded by the Monte Carlo sampling. The relative and absolute error tolerances were set to 0.001 and the calculations were performed using a 40 by 40 lattice. The algorithm converged after 5-7 iterations depending on the temperature. The condition number for the Jacobian matrix was about $10^4$ for all temperatures.

We now turn to the main focus of the present paper, which is the computation of marginal densities for known exponential densities. Suppose that we begin with an Ising spin square lattice of size $L.$ We split the vector of variables $\fx$ as $\fx=(\hat{\fx},\tilde{\fx}).$ The vector $\hat{\fx}$ contains the variables whose marginal density we want to compute and $\tilde{\fx}$ is the vector of variables we wish to eliminate.  We use the same 8 potential functions as above. The exponential density for the vector $\fx$ is given by $\exp(-H(\alpha,\fx)),$ where the Hamiltonian $H(\alpha,\fx)=\sum_{k=1}^{8} \alpha_k \psi_k (\fx),$ with $\alpha_1=2/T$ and $\alpha_l=0$ for $l=2,\ldots,8.$ We can apply the algorithm of Section \ref{kno} and obtain the density $\exp(-\hat{H}(\hat{\alpha},\hat{\fx}))/Z(\hat{\fa})$ for the vector $\hat{\fx}.$ The Hamiltonian is $\hat{H}(\hat{\alpha},\hat{\fx})=\sum_{k=1}^{8} \hat{\alpha}_k \hat{\psi}_k (\hat{\fx}),$ where the parameter vector $\hat{\alpha}$ contains the new values that were computed from the algorithm. Of course, once the density for $\hat{\fx}$ is known, we can repeat the procedure and find the density for a subset of the vector  $\hat{\fx}$ and so on.  We can make this recursive elimination of variables more systematic. Suppose that we start with a lattice of size $L\times L$ with $L$ even. We denote the vector of spins corresponding to this lattice as $\fx^{(0)}.$ We construct $2 \times 2$ blocks of variables and we represent each block by one variable, say the lower left-hand corner variable of the block. The vector of the new variables is denoted by $\fx^{(1)}$ and occupies a lattice of  size $L/2 \times L/2.$ There are different ways of assigning values to the new variables. Here we pick the "majority" rule. This means that if the sum of spins in a block is positive, the new variable is assigned the value of 1, if the sum of spins is negative the value -1 . In the case of a tie, the value of the  the new variable is taken to be the value of the spin on the lower left-hand corner. Then we apply the algorithm and we obtain the density for the new variables. As we have already mentioned, we are not sure beforehand, that the potential functions we choose are enough to describe all the couplings that may appear among the new variables. Thus the marginal density that we compute is usually only an approximation to the marginal density of the new variables. Once the density for the vector $\fx^{(1)}$ is computed, we can repeat the coarse-graining process and obtain the density for the vector $\fx^{(2)}$ which occupies a lattice of size $L/4 \times L/4$ and so on. The corresponding parameter vectors $\alpha^{(0)}, \alpha^{(1)}, \alpha^{(2)}, \ldots$ define a "flow" in the space of parameters. As is known \cite{binney}, the flow of parameter vectors contains useful information for a system near its critical point. In addition, we can store the parameter vectors for the different marginalization steps for future use.

We use the parameter vectors for the successive renormalization steps to locate the critical temperature of the Ising model. For temperatures $T < T_c$ the renormalization steps couple ever more distant spins while for $T > T_c$ the spins become successively decoupled. We can measure this successive coupling or decoupling by focusing on the quadratic terms (the terms $-\sum_J X_{J,1} X_{J,k+1}$) in the expansion of the Hamiltonian and calculating the "second moments" of their associated parameters
$$M_2^{(j)}= \sum_{k=1}^{l_q} d_k^2 \alpha^{(j)}_k$$
The superscript $j$ is used to denote the parameters at the $j$-th renormalization step and $M_2^{(0)}=\alpha_1^{(0)}$ is the second moment at the fine level ($j=0$). The parameters $d_k, \, k=1,\ldots,l_q$ denote the distance from $J$  of the spins in the group $k,$ while $l_q$ is the number of quadratic terms in the expansion of the Hamiltonian. For our experiments $l_q=5$ (see the definition of potential functions in (\ref{pots})). 
\begin{figure}
\centering
\epsfig{file=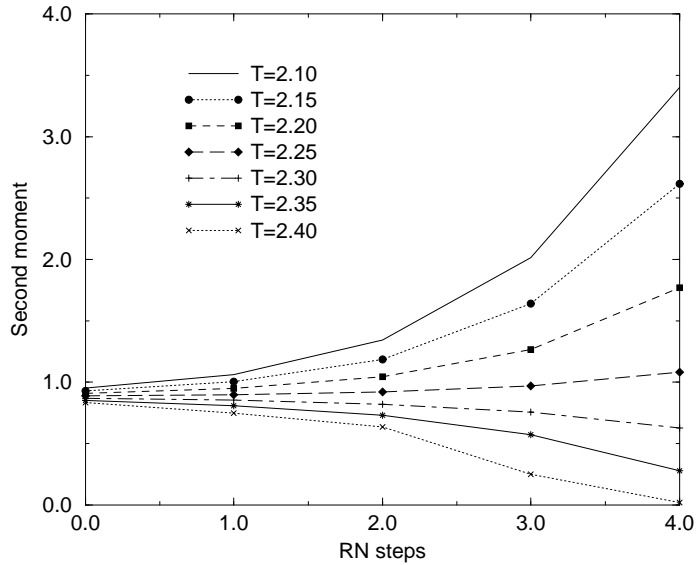,height=3in}
\caption{Second moment for successive renormalization steps for different temperatures.}
\label{fig:kno1}
\end{figure}
\begin{figure}
\centering
\epsfig{file=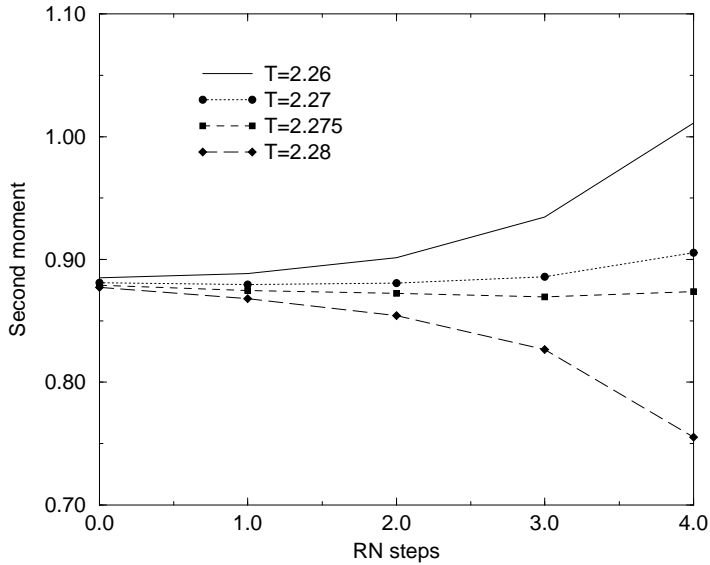,height=3in}
\caption{Second moment for successive renormalization steps for temperatures close to the critical.}
\label{fig:kno2}
\end{figure}
In Fig.(\ref{fig:kno1}) we show the evolution of $M_2^{(j)}$ for $j=0,1,\ldots,4$ and different temperatures. The fine level lattice is 80 by 80 and each renormalization step (based on the majority rule) decreases the number of spins by a factor of 4. For all the renormalization steps, the parameters were initialized at their values at the fine level. This means, that for $j=1,\ldots,4$ we started the algorithm with $\alpha_1^{(j)}=2/T$ and $\alpha_k^{(j)}=0, \, k=2,\ldots,8.$  As expected from the analytical results, below the critical temperature the second moment increases with successive renormalization steps, while for temperatures above the critical one, the second moment decreases with each renormalization step. The temperature for which the second moment remains constant (within the accuracy afforded by Monte Carlo sampling) should be the critical temperature. For our choice of potential functions and $10^5$ Monte Carlo steps per spin, the critical temperature is found to be $T_c \sim 2.275,$ an error of about $.3 \%.$ This can be seen in Fig.(\ref{fig:kno2}) where we have focused on a tighter temperature interval around the critical temperature. For $T=2.275,$ the second moment remains constant to within the accuracy of the Monte Carlo sampling. The relative and absolute tolerances were set to 0.001. The Levenberg-Marquardt algorithm converged after 5-12 iterations depending on the temperature and the renormalization step. As before, the condition number for the Jacobian matrix was about $10^4$ for all temperatures and renormalization steps.

\section*{Conclusions}
We have presented an algorithm for the estimation and renormalization of exponential densities. The algorithm is based on the maximization of the likelihood and results in a moment-matching problem. The matching problem is solved as an optimization problem with the Levenberg-Marquardt algorithm. For the case of renormalization, the moments for the reduced system are computed using Swendsen's renormalization method. We have exhibited the use of the algorithm by applying it to the Ising model, where it reproduces the analytical results with good accuracy. 

We hope that the algorithm can be useful in diverse areas where one needs to estimate an unknown density or marginalize a known density. Such areas range from applications to graphical models (the inference problem) to reductions of systems of differential equations. In the latter case, we are planning to apply the algorithm to estimate the density of solutions produced by finite-difference or spectral formulations of partial differential equations that do not have a natural candidate for a density. An analytical expression for the density in such situations can be helpful in the process of constructing reduced models \cite{chorin1}.

\section*{Acknowledgements}
I would like to thank Professor Alexandre Chorin and Professor Ole Hald for very helpful discussions and comments. This work was supported in part by the Applied Mathematical Sciences subprogram of the Office of Energy Research of the US Department of Energy under Contract DE-AC03-76-SF00098 and in part by the National Science Foundation under Grant DMS98-14631.

\bibliographystyle{elsart-num}

\bibliography{paper}

\end{document}